\documentclass[12pt,a4paper,reqno]{amsart}
\usepackage{amssymb,amsfonts}
\usepackage[english]{babel}
\usepackage{srcltx}
\usepackage{enumerate}
\advance\textwidth50mm \advance\hoffset-20mm
\advance\textheight30mm \advance\voffset-15mm
\newtheorem{theorem}{Theorem}

\newtheorem{corollary}{Corollary}
\theoremstyle{remark}

\begin{document}

\title{A Bernstein--Ganzburg limit theorem for best weighted approximation}
\author{D.~V.~Gorbachev}

\address{Lomonosov Moscow State University (Moscow)\\
Moscow Center of Fundamental and Applied Mathematics (Moscow)}

\email{dvgmail@mail.ru}

\thanks{The research was supported by the Russian Science Foundation (project
no.~23-71-30001) at Lomonosov Moscow State University.}

\date{}

\begin{abstract}
We prove a Bernstein--Ganzburg type limit relation
\[
\lim_{n\to\infty}
\Bigl(\frac{n}{\sigma}\Bigr)^{(2a+1)/p}E_{n,\sigma}(f)_{p,a,b}
=A_{\sigma}(f)_{p,a},
\]
where $E_{n,\sigma}(f)_{p,a,b}$ is the error of best approximation of
$f(nt/\sigma)$ by trigonometric polynomials of degree at most $n$ in
$L^{p}((-\pi,\pi],|2\sin(t/2)|^{2a}|\cos(t/2)|^{2b}\,dt)$, and
$A_{\sigma}(f)_{p,a}$ is the error of best approximation of $f$ by entire
functions of exponential type at most $\sigma$ in
$L^{p}(\mathbb{R},|x|^{2a}\,dx)$. For $a=b=0$, this result was obtained by
M.~I.~Ganzburg. The proof uses ideas from the Bernstein--Ganzburg limit
theorems and a localization method with the Fej\'er kernel from the proof of
the limit relation for Nikol'skii constants. As an application, using known
results for polynomial approximation, we compute the exact value of
$A_{\pi}(\mathbf{1}_{(-1,1)})_{1,a}$ for $a=0$ and $a=1/2$.
\end{abstract}

\maketitle

\section{Introduction and main result}

Let
\[
1\le p<\infty,\quad a\ge b\ge 0.
\]
We introduce the periodic Jacobi weight \cite{Vi22}
\[
w_{a,b}(t)= \Bigl|2\sin\frac{t}{2}\Bigr|^{2a}
\Bigl|\cos\frac{t}{2}\Bigr|^{2b},\quad t\in \mathbb{R}.
\]
The normalization by $2\sin(t/2)$ will be convenient in the estimates, since
$|2\sin(t/2)|\le |t|$.

Let
\[
L^{p}(w_{a,b})=L^{p}((-\pi,\pi],w_{a,b}(t)\,dt),\quad
L^{p}(|x|^{2a})=L^{p}(\mathbb{R},|x|^{2a}\,dx).
\]

By $\mathcal{T}_{n}$ we denote the space of $2\pi$-periodic trigonometric
polynomials of degree at most $n$, and by $\mathcal{E}_{\sigma}$ the class of
entire functions of exponential type at most~$\sigma$. Put
\[
\mathcal{E}_{\sigma,p,a}= \mathcal{E}_{\sigma}\cap L^{p}(|x|^{2a}).
\]
If $F\in\mathcal{E}_{\sigma,p,a}$, then the weighted Nikol'skii inequality (see,
for example, \cite{Go21}) implies that $F$ is bounded on $\mathbb{R}$. Therefore,
for all $z\in\mathbb{C}$ we have the estimate
\[
|F(z)|\le \|F\|_{L^{\infty}(\mathbb{R})} e^{\sigma|\!\operatorname{Im}z|}.
\]

We consider the problem of best approximation of a function
\[
f\in L^{p}(|x|^{2a}).
\]

For $n\in\mathbb{N}$ and $\sigma>0$, define the error of periodic approximation by
\[
E_{n,\sigma}(f)_{p,a,b}=\inf_{T\in\mathcal{T}_{n}}
\biggl(\int_{-\pi}^{\pi}\Bigl|f\Bigl(\frac{n}{\sigma}\,t\Bigr)-T(t)\Bigr|^{p}
w_{a,b}(t)\,dt\biggr)^{1/p}.
\]
The error on the line is defined by
\[
A_{\sigma}(f)_{p,a}= \inf_{F\in\mathcal{E}_{\sigma,p,a}}
\biggl(\int_{\mathbb{R}}|f(x)-F(x)|^{p} |x|^{2a}\,dx\biggr)^{1/p}.
\]

\begin{theorem}\label{thm-main}
Let $1\le p<\infty$, $a\ge b\ge 0$, and $\sigma>0$. For every function $f\in
L^{p}(|x|^{2a})$ one has
\[
\lim_{n\to\infty}
\Bigl(\frac{n}{\sigma}\Bigr)^{(2a+1)/p}E_{n,\sigma}(f)_{p,a,b}=A_{\sigma}(f)_{p,a}.
\]
\end{theorem}

This statement is a weighted periodic analogue of the Bernstein--Ganzburg limit
theorems for errors of best approximation by algebraic polynomials and by entire
functions of exponential type. In \cite{Ga21,Ga23}, such relations were developed
in a general form for errors of best approximation and sharp constants. However,
in our case only the unweighted case $a=b=0$ was covered. On the other hand, the
proof below is close to the method for limit estimates of Nikol'skii constants
\cite{Go21}: in both problems one uses local scaling of the Jacobi weight and
multiplication by the Fej\'er kernel.

The key point of the passage to the limit is that, for fixed $x\in\mathbb{R}$,
\[
\Bigl(\frac{n}{\sigma}\Bigr)^{2a} w_{a,b}\Bigl(\frac{\sigma x}{n}\Bigr)\to
|x|^{2a},\quad n\to \infty,
\]
and also
\[
\Bigl(\frac{n}{\sigma}\Bigr)^{2a} w_{a,b}\Bigl(\frac{\sigma x}{n}\Bigr)\le
|x|^{2a}.
\]
Hence, in particular,
\begin{equation}\label{w}
\Bigl(\frac{n}{\sigma}\Bigr)^{(2a+1)/p}
\Bigl\|f\Bigl(\frac{n}{\sigma}\,{\,\cdot\,}\Bigr)\Bigr\|_{L^{p}(w_{a,b})}=
\biggl(\int_{-\pi n/\sigma}^{\pi n/\sigma}|f(x)|^{p}
\Bigl(\frac{n}{\sigma}\Bigr)^{2a}w_{a,b}\Bigl(\frac{\sigma
x}{n}\Bigr)\,dx\biggr)^{1/p}\le \|f\|_{L^{p}(|x|^{2a})}.
\end{equation}

In what follows, for brevity we write
\[
E_{n,\sigma}(f)= E_{n,\sigma}(f)_{p,a,b},\quad
A_{\sigma}(f)=A_{\sigma}(f)_{p,a}.
\]

\section{Lower estimate}

In this section we prove
\begin{equation}\label{A-sigma-E}
A_{\sigma}(f)\le \liminf_{n\to\infty}
\Bigl(\frac{n}{\sigma}\Bigr)^{(2a+1)/p}E_{n,\sigma}(f).
\end{equation}
It is convenient first to consider the two-parameter error
\[
E_{n,N,\sigma}(f)= \inf_{T\in\mathcal{T}_{n}}
\biggl(\int_{-\pi}^{\pi}\Bigl|f\Bigl(\frac{N}{\sigma}\,t\Bigr)-T(t)\Bigr|^{p}
w_{a,b}(t)\,dt\biggr)^{1/p}.
\]
At the end we return to the case $N=n$.

We introduce the integral Fej\'er kernel
\begin{equation}\label{varphi1}
\varphi(t)=\Bigl(\frac{\sin(t/2)}{t/2}\Bigr)^{2},\quad \varphi(0)=1,
\end{equation}
which is an entire function of exponential type~$1$. We shall use the fact that
\begin{equation}\label{varphi2}
0\le 1-\varphi(t)\le \min\,\Bigl\{\frac{t^{2}}{12},1\Bigr\},\quad t\in
\mathbb{R}.
\end{equation}

Choose $m=m_{n}\in \mathbb{N}$ such that
\[
m_{n}\to\infty, \quad m_{n}=o(n),\quad n\to \infty,\quad 2mp\ge2a+2,
\]
for example $m_{n}=\lceil\!\sqrt{n}+c\,\rceil$ with sufficiently large $c$.
Fix $q\in \mathbb{N}$ satisfying
\[
2qp\ge 2b.
\]
Put
\[
N=N_{n}=n+m_{n}+q.
\]
Then $N_{n}/n\to 1$.

Consider the smoothing function
\[
\eta_{m}(t)=\varphi(t)^{m}\Bigl(\cos\frac{t}{2}\Bigr)^{2q}.
\]
It is an entire function of exponential type $m+q$ such that, for $t\in \mathbb{R}$,
\begin{equation}\label{eta1}
0\le \eta_{m}(t)\le 1.
\end{equation}
Taking into account \eqref{varphi2} and the elementary inequality
\[
1-u^{r}\le r(1-u),\quad u\in [0,1],\quad r\ge 1,
\]
we get
\begin{equation}\label{eta2}
1-\eta_{m}(t)\le 1-\varphi(t)^{m}+ 1-\Bigl(\cos\frac{t}{2}\Bigr)^{2q}\le
m(1-\varphi(t)) +q\Bigl(\sin\frac{t}{2}\Bigr)^{2}\le
\Bigl(\frac{m}{12}+\frac{q}{4}\Bigr)t^{2}.
\end{equation}

For an arbitrary trigonometric polynomial $T\in\mathcal{T}_{n}$ define an
entire function of exponential type~by
\[
G(x)=G_{T,m,N,\sigma}(x)= T\Bigl(\frac{\sigma
x}{N}\Bigr)\eta_{m}\Bigl(\frac{\sigma x}{N}\Bigr).
\]
The function $T(\sigma x/N)$ has type at most $n\sigma/N$, while
$\eta_{m}(\sigma x/N)$ has type $(m+q)\sigma/N$. Hence $G$ has type at most
$(n+m+q)\sigma/N=\sigma$. Moreover,
\[
|G(x)|^{p}\,|x|^{2a}=O(|x|^{-2mp+2a}),\quad |x|\to\infty.
\]
Therefore, since $2mp\ge2a+2$,
\[
G\in\mathcal{E}_{\sigma,p,a}.
\]

For all $t\in \mathbb{R}$, by $2mp\ge2a+2$ and $2qp\ge2b$, we have
\begin{align*}
\eta_{m}(t)^{p}|t|^{2a}&= \Bigl|\frac{\sin(t/2)}{t/2}\Bigr|^{2mp}
\Bigl|\cos\frac{t}{2}\Bigr|^{2qp}|t|^{2a}\\ &=
\Bigl|2\sin\frac{t}{2}\Bigr|^{2a} \Bigl|\frac{\sin(t/2)}{t/2}\Bigr|^{2mp-2a}
\Bigl|\cos\frac{t}{2}\Bigr|^{2qp}\\ &\le
\Bigl|2\sin\frac{t}{2}\Bigr|^{2a}\Bigl|\cos\frac{t}{2}\Bigr|^{2b}= w_{a,b}(t).
\end{align*}
Thus, making the change of variable $t=\sigma x/N$, we obtain the inequality
\begin{align}
\int_{-\pi N/\sigma}^{\pi N/\sigma} \Bigl|f(x)\eta_{m}\Bigl(\frac{\sigma
x}{N}\Bigr)-G(x)\Bigr|^{p} |x|^{2a}\,dx&= \Bigl(\frac{N}{\sigma}\Bigr)^{2a+1}
\int_{-\pi}^{\pi} \Bigl|f\Bigl(\frac{N}{\sigma}\,t\Bigr)-T(t)\Bigr|^{p}
\eta_{m}(t)^{p} |t|^{2a}\,dt \notag\\ &\le \Bigl(\frac{N}{\sigma}\Bigr)^{2a+1}
\int_{-\pi}^{\pi} \Bigl|f\Bigl(\frac{N}{\sigma}\,t\Bigr)-T(t)\Bigr|^{p}
w_{a,b}(t)\,dt. \label{ineq1}
\end{align}

Next, for $t\in [-\pi,\pi]$ and $k\in \mathbb{Z}\setminus \{0\}$, we have
\begin{align*}
\eta_{m}(t+2\pi k)^{p} |t+2\pi k|^{2a}&= \Bigl|\frac{2\sin(t/2)}{t+2\pi
k}\Bigr|^{2mp} \Bigl|\cos\frac{t}{2}\Bigr|^{2qp} |t+2\pi k|^{2a}\\ &=
\Bigl|\frac{2\sin(t/2)}{t+2\pi k}\Bigr|^{2mp-2a}
\Bigl|2\sin\frac{t}{2}\Bigr|^{2a} \Bigl|\cos\frac{t}{2}\Bigr|^{2qp}\\ &\le
\Bigl|\frac{2}{t+2\pi k}\Bigr|^{2mp-2a}\Bigl|2\sin\frac{t}{2}\Bigr|^{2a}
\Bigl|\cos\frac{t}{2}\Bigr|^{2b}\\ &\le
\Bigl(\frac{2}{\pi(2|k|-1)}\Bigr)^{2mp-2a}w_{a,b}(t).
\end{align*}
Hence we find
\begin{align}
\int_{|x|>\pi N/\sigma} |G(x)|^{p} |x|^{2a}\,dx&=
\Bigl(\frac{N}{\sigma}\Bigr)^{2a+1} \int_{|t|>\pi} |T(t)|^{p} \eta_{m}(t)^{p}
|t|^{2a}\,dt \notag\\ &= \Bigl(\frac{N}{\sigma}\Bigr)^{2a+1} \sum_{k\ne
0}\int_{[-\pi,\pi]+2\pi k} |T(t)|^{p} \eta_{m}(t)^{p} |t|^{2a}\,dt \notag\\ &=
\Bigl(\frac{N}{\sigma}\Bigr)^{2a+1} \int_{-\pi}^{\pi}|T(t)|^{p} \sum_{k\ne0}
\eta_{m}(t+2\pi k)^{p} |t+2\pi k|^{2a}\,dt \notag\\ &\quad \le \delta_{m}^{p}
\Bigl(\frac{N}{\sigma}\Bigr)^{2a+1} \int_{-\pi}^{\pi}|T(t)|^{p} w_{a,b}(t)\,dt,
\label{ineq2}
\end{align}
where
\[
\delta_{m}^{p}= \sum_{k\ne 0}
\Bigl(\frac{2}{\pi(2|k|-1)}\Bigr)^{2mp-2a}\le
\Bigl(\frac{2}{\pi}\Bigr)^{2mp-2a} \sum_{k\ne 0}\frac{1}{(2|k|-1)^{2}}
=
\frac{\pi^{2}}{4}\Bigl(\frac{2}{\pi}\Bigr)^{2mp-2a}.
\]
We have
\[
\delta_{m}\to 0,\quad m\to\infty.
\]

We now prove \eqref{A-sigma-E}. Put
\begin{equation}\label{L}
L= \liminf_{n\to\infty} \Bigl(\frac{N_{n}}{\sigma}\Bigr)^{(2a+1)/p}
E_{n,N_{n},\sigma}(f).
\end{equation}
If $L=+\infty$, there is nothing to prove. Otherwise, choose a subsequence on
which $L$ is attained. Without changing the notation, we shall assume that
\[
\Bigl(\frac{N_{n}}{\sigma}\Bigr)^{(2a+1)/p} E_{n,N_{n},\sigma}(f)\to L,\quad
n\to \infty.
\]

For each $n$ from this subsequence, choose $T_{n}\in\mathcal{T}_{n}$ such that
\[
\Bigl(\frac{N_{n}}{\sigma}\Bigr)^{(2a+1)/p}
\Bigl\|f\Bigl(\frac{N_{n}}{\sigma}\,{\,\cdot\,}\Bigr)-
T_{n}\Bigr\|_{L^{p}(w_{a,b})}\le L+\varepsilon_{n},
\]
where $\varepsilon_{n}>0$ and $\varepsilon_{n}\to0$. Now put
\[
G_{n,\sigma}=G_{T_{n},m_{n},N_{n},\sigma}.
\]

Since $G_{n,\sigma}\in\mathcal{E}_{\sigma,p,a}$, we have
\[
A_{\sigma}(f)\le \|f-G_{n,\sigma}\|_{L^{p}(|x|^{2a})}.
\]
Write
\[
f(x)-G_{n,\sigma}(x)= f(x)\eta_{m_{n}}\Bigl(\frac{\sigma x}{N_{n}}\Bigr)-
G_{n,\sigma}(x)+ f(x)\Bigl(1-\eta_{m_{n}}\Bigl(\frac{\sigma x}{N_{n}}\Bigr)\Bigr).
\]

We use \eqref{eta1} and \eqref{eta2}. Then, for every fixed
$x\in \mathbb{R}$, in view of $m_{n}=o(n)$ and $N_{n}/n\to1$,
\[
\eta_{m_{n}}\Bigl(\frac{\sigma x}{N_{n}}\Bigr)\to 1,
\]
and therefore, by the Lebesgue dominated convergence theorem,
\[
\varepsilon_{n,1}(f)= \Bigl\|f\Bigl(1-
\eta_{m_{n}}\Bigl(\frac{\sigma}{N_{n}}\,{\,\cdot\,}\Bigr)\Bigr)\Bigr\|_{L^{p}(|x|^{2a})}\to
0.
\]

From \eqref{eta1}, \eqref{ineq1}, and \eqref{ineq2} it follows that
\begin{align*}
\Bigl\|f\eta_{m_{n}}\Bigl(\frac{\sigma}{N_{n}}\,{\,\cdot\,}\Bigr)-
G_{n,\sigma}\Bigr\|_{L^{p}(|x|^{2a})}&\le
\Bigl\|\Bigl(f\eta_{m_{n}}\Bigl(\frac{\sigma}{N_{n}}\,{\,\cdot\,}\Bigr)-
G_{n,\sigma}\Bigr)\mathbf{1}_{\{|x|\le \pi
N_{n}/\sigma\}}\Bigr\|_{L^{p}(|x|^{2a})}\\ &\quad
{}+\Bigl\|\Bigl(f\eta_{m_{n}}\Bigl(\frac{\sigma}{N_{n}}\,{\,\cdot\,}\Bigr)-
G_{n,\sigma}\Bigr)\mathbf{1}_{\{|x|>\pi
N_{n}/\sigma\}}\Bigr\|_{L^{p}(|x|^{2a})}\\ &\le
\Bigl(\frac{N_{n}}{\sigma}\Bigr)^{(2a+1)/p}
\Bigl\|f\Bigl(\frac{N_{n}}{\sigma}\,{\,\cdot\,}\Bigr)-T_{n}\Bigr\|_{L^{p}(w_{a,b})}\\
&\quad {}+\|f\mathbf{1}_{\{|x|>\pi N_{n}/\sigma\}}\|_{L^{p}(|x|^{2a})}\\ &\quad
{}+ \delta_{m_{n}} \Bigl(\frac{N_{n}}{\sigma}\Bigr)^{(2a+1)/p}
\|T_{n}\|_{L^{p}(w_{a,b})}.
\end{align*}
Here
\[
\varepsilon_{n,2}(f)=\|f\mathbf{1}_{\{|x|>\pi N_{n}/\sigma\}}\|_{L^{p}(|x|^{2a})}\to 0,\quad n\to
\infty,
\]
and by \eqref{w}
\begin{align*}
\Bigl(\frac{N_{n}}{\sigma}\Bigr)^{(2a+1)/p} \|T_{n}\|_{L^{p}(w_{a,b})}&\le
\Bigl(\frac{N_{n}}{\sigma}\Bigr)^{(2a+1)/p}
\Bigl\|f\Bigl(\frac{N_{n}}{\sigma}\,{\,\cdot\,}\Bigr)-T_{n}\Bigr\|_{L^{p}(w_{a,b})}\\
&\quad {}+ \Bigl(\frac{N_{n}}{\sigma}\Bigr)^{(2a+1)/p}
\Bigl\|f\Bigl(\frac{N_{n}}{\sigma}\,{\,\cdot\,}\Bigr)\Bigr\|_{L^{p}(w_{a,b})}\\
&\le L+\varepsilon_{n}+\|f\|_{L^{p}(|x|^{2a})}.
\end{align*}

Thus,
\[
\|f-G_{n,\sigma}\|_{L^{p}(|x|^{2a})}\le
\varepsilon_{n,1}(f)+L+\varepsilon_{n}+\varepsilon_{n,2}(f)+
\delta_{m_{n}}(L+\varepsilon_{n}+\|f\|_{L^{p}(|x|^{2a})}),
\]
whence
\[
A_{\sigma}(f)\le \limsup_{n\to\infty} \|f-G_{n,\sigma}\|_{L^{p}(|x|^{2a})} \le
L.
\]

It remains to replace $N_{n}$ by $n$. We use the standard fact that, if $f\in
L^{p}(|x|^{2a})$, then
\[
\lim_{\lambda\to1,\ \lambda>0}
\|f(\lambda{\,\cdot\,})-f\|_{L^{p}(|x|^{2a})}=0.
\]

Put $\lambda_{n}=N_{n}/n\to 1$. Then, taking into account \eqref{w},
\begin{align*}
\Bigl(\frac{n}{\sigma}\Bigr)^{(2a+1)/p}
|E_{n,N_{n},\sigma}(f)-E_{n,\sigma}(f)|&\le
\Bigl(\frac{n}{\sigma}\Bigr)^{(2a+1)/p}\Bigl\|f\Bigl(\frac{N_{n}}{\sigma}\,t\Bigr)-
f\Bigl(\frac{n}{\sigma}\,t\Bigr)\Bigr\|_{L^{p}(w_{a,b})}\\ &\le
\|f(\lambda_{n}\,{\,\cdot\,})- f\|_{L^{p}(|x|^{2a})}\to 0.
\end{align*}
Hence, together with \eqref{L}, we get
\[
\liminf_{n\to\infty} \Bigl(\frac{n}{\sigma}\Bigr)^{(2a+1)/p}E_{n,\sigma}(f)=
\liminf_{n\to\infty}
\Bigl(\frac{N_{n}}{\sigma}\Bigr)^{(2a+1)/p}E_{n,N_{n},\sigma}(f)=L,
\]
which implies the desired inequality \eqref{A-sigma-E}.

\section{Upper estimate}

We now prove
\[
\limsup_{n\to\infty} \Bigl(\frac{n}{\sigma}\Bigr)^{(2a+1)/p}E_{n,\sigma}(f) \le
A_{\sigma}(f).
\]
Fix $\varepsilon>0$. Choose $F\in\mathcal{E}_{\sigma,p,a}$ such that
\[
\|f-F\|_{L^{p}(|x|^{2a})}\le A_{\sigma}(f)+\varepsilon.
\]

We construct a trigonometric polynomial from $F$. Recall that the Fej\'er kernel
$\varphi$ is defined by~\eqref{varphi1}.~Put
\[
T_{n}(t)= \sum_{k\in\mathbb{Z}} F\Bigl(\frac{n}{\sigma}\,(t+2\pi k)\Bigr)
\varphi(t+2\pi k).
\]
The function $F(nt/\sigma)\varphi(t)$ has type at most $n+1$ and is
$O(|t|^{-2})$ as $|t|\to \infty$. Therefore its Fourier transform is continuous
and, by the Paley--Wiener theorem, vanishes on $\mathbb{R}\setminus
(-(n+1),n+1)$, and hence also at all integers $k$ such that $|k|>n$.
Consequently, by the Poisson summation formula, $T_{n}\in\mathcal{T}_{n}$.

We have
\[
f\Bigl(\frac{n}{\sigma}\,t\Bigr)-T_{n}(t)=
\Bigl(f\Bigl(\frac{n}{\sigma}\,t\Bigr)-F\Bigl(\frac{n}{\sigma}\,t\Bigr)\Bigr) +
\Bigl(F\Bigl(\frac{n}{\sigma}\,t\Bigr)-T_{n}(t)\Bigr).
\]
For the first summand,
\begin{align*}
\Bigl(\frac{n}{\sigma}\Bigr)^{2a+1} \int_{-\pi}^{\pi}
\Bigl|f\Bigl(\frac{n}{\sigma}\,t\Bigr)-F\Bigl(\frac{n}{\sigma}\,t\Bigr)\Bigr|^{p}
w_{a,b}(t)\,dt&= \int_{-\pi n/\sigma}^{\pi n/\sigma}|f(x)-F(x)|^{p}
\Bigl(\frac{n}{\sigma}\Bigr)^{2a} w_{a,b}\Bigl(\frac{\sigma x}{n}\Bigr)\,dx\\
&\le \int_{\mathbb{R}}|f(x)-F(x)|^{p} |x|^{2a}\,dx.
\end{align*}

For the second summand we use the decomposition
\[
F\Bigl(\frac{n}{\sigma}\,t\Bigr)-T_{n}(t)=C_{n}(t)-R_{n}(t),
\]
where
\[
C_{n}(t)=F\Bigl(\frac{n}{\sigma}\,t\Bigr)(1-\varphi(t)),
\]
and
\[
R_{n}(t)= \sum_{k\ne0} F\Bigl(\frac{n}{\sigma}\,(t+2\pi k)\Bigr) \varphi(t+2\pi
k).
\]

First estimate the central term. After the change of variable $x=nt/\sigma$ we get
\[
J_{n}=\Bigl(\frac{n}{\sigma}\Bigr)^{2a+1} \|C_{n}\|_{L^{p}(w_{a,b})}^{p}=
\int_{-\pi n/\sigma}^{\pi n/\sigma}|F(x)|^{p} \Bigl|1-\varphi\Bigl(\frac{\sigma
x}{n}\Bigr)\Bigr|^{p} \Bigl(\frac{n}{\sigma}\Bigr)^{2a} w_{a,b}\Bigl(\frac{\sigma
x}{n}\Bigr)\,dx.
\]
By \eqref{varphi2},
\[
J_{n}\le \int_{\mathbb{R}}|F(x)|^{p} |x|^{2a}
\min\,\Bigl\{\Bigl(\frac{\sigma^{2}x^{2}}{12n^{2}}\Bigr)^{p},1\Bigr\}\,dx.
\]
Hence $J_{n}\to0$. Equivalently,
\[
\Bigl(\frac{n}{\sigma}\Bigr)^{(2a+1)/p}
\|C_{n}\|_{L^{p}(w_{a,b})}\to0.
\]

Now estimate the tail. Since
\[
\sum_{k\in\mathbb{Z}}\varphi(t+2\pi k)=1,\quad
t\in \mathbb{R},
\]
by H\"older's inequality we have
\begin{align*}
|R_{n}(t)|^p&= \biggl|\sum_{k\ne0} F\Bigl(\frac{n}{\sigma}\,(t+2\pi
k)\Bigr)\varphi(t+2\pi k)\biggr|^p \\ &\le \biggl(\sum_{k\ne0}\varphi(t+2\pi
k)\biggr)^{p-1} \sum_{k\ne0} \Bigl|F\Bigl(\frac{n}{\sigma}\,(t+2\pi
k)\Bigr)\Bigr|^p \varphi(t+2\pi k) \\ &\le \sum_{k\ne0}
\Bigl|F\Bigl(\frac{n}{\sigma}\,(t+2\pi k)\Bigr)\Bigr|^p \varphi(t+2\pi k).
\end{align*}
Integrating and changing the variable $x=n(t+2\pi k)/\sigma$, we get
\begin{align*}
\Bigl(\frac{n}{\sigma}\Bigr)^{2a+1} \|R_{n}\|_{L^{p}(w_{a,b})}^{p}&\le
\int_{|x|>\pi n/\sigma}|F(x)|^{p} \varphi\Bigl(\frac{\sigma x}{n}\Bigr)
\Bigl(\frac{n}{\sigma}\Bigr)^{2a} w_{a,b}\Bigl(\frac{\sigma x}{n}\Bigr)\,dx\\
&\le \int_{|x|>\pi n/\sigma}|F(x)|^{p}|x|^{2a}\,dx \to0.
\end{align*}
Therefore,
\[
\Bigl(\frac{n}{\sigma}\Bigr)^{(2a+1)/p}
\Bigl\|F\Bigl(\frac{n}{\sigma}\,t\Bigr)-T_{n}(t)\Bigr\|_{L^{p}(w_{a,b})}
\to0.
\]

As a result,
\[
\limsup_{n\to\infty} \Bigl(\frac{n}{\sigma}\Bigr)^{(2a+1)/p}E_{n,\sigma}(f)
\le
\|f-F\|_{L^{p}(|x|^{2a})} \le A_{\sigma}(f)+\varepsilon.
\]
Since $\varepsilon>0$ is arbitrary, the upper estimate is proved. Together with
the lower estimate, this completes the proof of the main theorem.

We note the following assertion useful for applications, which follows from the
proof of the lower estimate. The notation is also taken from there.

\begin{corollary}\label{cor-main}
Let $N_n>0$, $N_n'>0$, $N_n/n\to1$, and $N_n'/n\to1$. Then, under the assumptions
of Theorem~\ref{thm-main},
\[
\lim_{n\to\infty} \Bigl(\frac{N_n'}{\sigma}\Bigr)^{(2a+1)/p}
E_{n,N_n,\sigma}(f)_{p,a,b}= A_{\sigma}(f)_{p,a}.
\]
\end{corollary}

\section{An application to approximation of individual functions}

Theorem~\ref{thm-main} gives a limit relation between the errors of best
weighted approximation by trigonometric polynomials and by entire functions of
exponential type. If we know the answer in one of the problems, then we can
immediately derive a corresponding result for the other one. Some examples are
given in \cite{Ga21}.

Here, as a consequence of Theorem~\ref{thm-main}, we obtain exact values of
$A_{\pi}(\chi_{1})_{1,a}$ for the characteristic function of the interval
\[
\chi_{r}(x)=\mathbf{1}_{(-r,r)}(x),\quad r>0,
\]
in the cases $a=b=0$ and $a=b=1/2$.

\subsection*{The case $a=b=0$}
We have
\[
A_{\pi}(\chi_{1})_{1,0}=\inf_{F\in\mathcal{E}_{\pi,1,0}}
\int_{\mathbb{R}}|\chi_{1}(x)-F(x)|\,dx.
\]

Denote $I_{n}(h)=E_{n-1}(\chi_{h})_{1,0,0}$. By Corollary~\ref{cor-main},
\[
A_{\pi}(\chi_{1})_{1,0}=\lim_{n\to\infty}\frac{n}{\pi}\,I_{n}\Bigl(\frac{\pi}{n}\Bigr).
\]
In \cite{Ba08} it was proved that
\[
\lim_{n\to\infty}nI_{n}\Bigl(\frac{\pi}{n}\Bigr)=2\pi-4v_0,
\]
where $v_0$ is the unique root of the equation
\[
\cos v=\frac{2\pi v}{v^2+\pi^2},\quad 0<v<\frac{\pi}{2}.
\]
Hence
\[
A_{\pi}(\chi_{1})_{1,0}=2-\frac{4v_0}{\pi}.
\]
Numerically,
\[
v_0=0.97116830789\ldots,\quad A_{\pi}(\chi_{1})_{1,0}=0.7634701058\ldots.
\]

\subsection*{The case $a=b=1/2$}
We have $w_{1/2,1/2}(t)=|\!\sin t|$ and
\[
A_{\pi}(\chi_{1})_{1,1/2}=\inf_{F\in\mathcal{E}_{\pi,1,1/2}} \int_{\mathbb{R}}
|\chi_{1}(x)-F(x)|\,|x|\,dx.
\]

Denote $J_{n}(h)=E_{n}(\chi_{h})_{1,1/2,1/2}$. Since the function $\chi_{h}$
and the weight $|\!\sin t|$ are even, the extremal polynomial may be assumed to
be even. Therefore, after the change of variable $u=\cos t$, we get
\[
J_{n}(h)=2\inf_{P\in\mathcal P_{n}}\int_{-1}^{1}|\zeta_{\cos h}(u)-P(u)|\,du,
\]
where
\[
\zeta_{\alpha}(u)=\mathbf{1}_{[\alpha,1]}(u),\quad -1<\alpha<1,
\]
and $P\in \mathcal P_{n}$ is an algebraic polynomial of degree at most $n$.

In \cite{De10} it was proved that, for
\[
\alpha=\cos\frac{\pi}{n+2},
\]
one has
\[
\inf_{P\in\mathcal P_{n}}\int_{-1}^{1}|\zeta_{\alpha}(u)-P(u)|\,du
=1-\alpha.
\]
Consequently,
\[
J_{n}\Bigl(\frac{\pi}{n+2}\Bigr)
=2\Bigl(1-\cos\frac{\pi}{n+2}\Bigr).
\]

Applying Corollary~\ref{cor-main}, we obtain
\[
A_{\pi}(\chi_{1})_{1,1/2}= \lim_{n\to\infty} \Bigl(\frac{n+2}{\pi}\Bigr)^{2}
J_{n}\Bigl(\frac{\pi}{n+2}\Bigr)= \lim_{n\to\infty}
\Bigl(\frac{n+2}{\pi}\Bigr)^{2} 2\Bigl(1-\cos\frac{\pi}{n+2}\Bigr)=1.
\]

In conclusion, let us say that it would be interesting to compute the value
$A_{\pi}(\chi_{1})_{1,a}$ for all $a\ge 0$. Such results are known, for
example, for one-sided approximations of the characteristic function~$\chi_{1}$
(see \cite{Ba17,Bu15,Ca17}). However, in our case the problem seems especially
difficult.

\end{document}